\documentclass[12pt,letterpaper]{amsart}    
\usepackage{amssymb,ifthen,color}
\usepackage{amsmath,amsthm}
\usepackage{latexsym}
\usepackage{enumerate}
\usepackage{amsfonts}
\usepackage{url}
\usepackage{graphicx}

\begin{document}

\title{Polish  Mathematicians and Mathematics in World War I}

\author[S. Domoradzki]{Stanis\l aw Domoradzki}
\address{Faculty of Mathematics and Natural Sciences, University of
Rzesz\'ow,  Ul. Prof. S. Pigonia  1, 35-959, Rzesz\'ow, Poland}
\email{domoradz@ur.edu.pl}
\author[M. Stawiska]{Ma{\l}gorzata Stawiska}
\address{Mathematical Reviews, 416 Fourth St., Ann Arbor, MI 48103, USA}
\email{stawiska@umich.edu}

\date{\today}                                           

\maketitle


 \tableofcontents

\setcounter{tocdepth}{3}

\begin{abstract}
In this article we  present  diverse experiences of Polish mathematicians (in a broad sense) who during World War I  fought for freedom of their homeland or  conducted their research  and teaching in difficult wartime circumstances. We first focus on  those affiliated with  Polish institutions of higher education: the existing Universities in Lw\'ow in Krak\'ow and the Lw\'ow Polytechnics  (Austro-Hungarian empire) as well as the reactivated University of Warsaw and the new Warsaw Polytechnics (the Polish Kingdom, formerly in the Russian empire). Then we consider  the situations of Polish mathematicians in the Russian empire and other countries. We discuss  not only individual fates, but also organizational efforts of many kinds  (teaching at the academic level outside traditional institutions-- in Society for Scientific Courses in Warsaw and in Polish University College in Kiev; scientific societies in Krak\'ow, Lw\'ow, Moscow and Kiev; publishing activities) in order to illustrate the formation of modern Polish mathematical community. 
\end{abstract}

\section{Introduction} When thinking about Polish mathematicians  during World War I one needs to keep in mind that from  1795 until 1918 there was no Poland on the political map of Europe. The territory of the former Crown of the Kingdom of Poland and the Grand Duchy of Lithuania (also known as the Commonwealth of Poland) was divided among three neighboring powers: Austria (later Austro-Hungary), Prussia (later part of the unified Germany) and Russia.   In 1866  the Austro-Hungarian partition (known as the province of Galicja)  was granted autonomy and Polish became the official language.  Polish-language academic centers (the only ones in the world) were located in
 Krak\'ow (Jagiellonian University) and Lw\'ow (University and Polytechnic). 
The Russian partition (the Polish Kingdom) saw periods of struggle for national independence (November Uprising 1830-1831 and January Uprising 1863-64) followed by repressions.  All public high schools conducted instruction in Russian and the   Main School in Warsaw gave way to the Russian-language Imperial University and Polytechnic. The use of Polish in official communication was forbidden.  The administration and courts were staffed by Russians.  Educated Poles had to look for career opportunities elsewhere, and many of them established themselves in other places in the Russian Empire. But the end of the 19th century in the Kingdom brought an amazing development of unofficial education in Polish at all levels, including academic. In the Prussian partition Poles were subject to national and religious dscrimination; the efforts at Germanization intensified after the unification of Germany in 1871. There were no academic institutions there.\footnote{Wroc\l aw (Breslau) had an  university, but the city was not within Polish borders in the 18th century.} So whom do we consider a ``Polish" mathematician? As others pointed  out before (\cite{dudamat}, \cite{obczyzna}), there is no good   answer to this question. This article concerns mostly  men and women born in the Polish Kingdom or Galicja and educated there (at least up to the high-school level), speaking Polish as one of their primary languages, including a few who were ethnically Jewish. But we also talk about some individuals born in the Russian empire to Polish-speaking families. We organized the article around places where Polish mathematicians found themselves  during the war and we discuss each person  in connection with the place of his or her  main wartime activity. We occasionally extend the discussion beyond November 11, 1918, the date of armistice and  proclamation of independent Poland. The borders of the reborn state were not guaranteed, so Poland had to defend itself against the competing interests of its neighbors: Germany, Russia, Lithuania, Czechoslovakia, Soviet republics   and multiple short-lived Ukrainian states. The fighting continued until 1921 (the Peace of Riga with Soviet Russia and Soviet Ukraine as well as the Third Silesian Uprising against Germany).\footnote{Stanis\l aw Saks fought in Silesian Uprisings; Antoni \L omnicki and Kazimierz Bartel fought in the Defense of Lw\'ow; Stanis\l aw Le\'sniewski, Stefan Mazurkiewicz and Wac\l aw Sierpi\'nski worked in a cryptography  group during the Polish-Soviet war.} \\

When the World War I broke out on July 28, 1914, it pitched Germany and Austro-Hungary against Russia (backed by France and Great Britain). Poles could  either  remain passive or  take sides in the conflict. Many were drafted into respective armies.\footnote{The total  number of men mobilized by Germany during WWI was 13.2 million, or 41.4\% of the male population. In Russia, the number was 13 million men, or 7.4\% of the male population, plus 5,000-6,000 women (\cite{beckett}).} Some did not believe in fighting a war serving the interests of the occupying powers; a war viewed as fratricidal, given the Polish presence in the enemy states. On the other hand, many saw an opportunity to fight for the independent Polish state aligning themselves with the side of either Central Powers or of Triple Entente. This was particularly true in Galicja, where young men massively volunteered into Polish Legions.\footnote{Only those not subject to draft into the regular Austro-Hungarian army were allowed to enlist in the Legions.  Some Poles from the Russian partiton studying in Galicja  did so (\cite{kutrzeba}).} Mathematicians were among those who  were drafted or volunteered into military, who  served in the trenches or in non-combat units, experiencing  wounds, gas poisoning, prisoners' camps or internment. They were also affected by compulsory evacuations\footnote{The Russians evacuated about 130 industrial enterprises and  200 educational institutions-- their personnel, equipment etc.--from the Polish  Kingdom and part of Eastern Galicja, over   600 thousand people in total.}, travel restrictions, food and raw materials shortages. None lost their life as a result of war operations, but two outstanding scholars (Marian Smoluchowski and Zygmunt Janiszewski) died of war-related epidemics.  Although the fates  we describe  were typical in many ways,  presenting the realities  of war from  personal perspectives  is not the only goal of our article. Rather, we want to focus on the impact of war on activities considered normal for mathematicians and other scientists: studying, research, teaching, academic administration,  publishing and professional organizations. We found it convenient to adopt a somewhat broad definition of a ``mathematician".  Besides scholars known for their  outstanding results in mathematics (prior to the war or afterwards), such as Banach and Sierpi\'nski,  we introduce many individuals who made lesser contributions to mathematical knowledge, in particular those (Izabela Abramowicz, Zygmunt Chwia\l kowski, Adam Patryn)  who did not continue their research  after the war.  Finally, we mention some physicists, astronomers, engineers and philosophers who in the circumstances of the war engaged in  teaching mathematics at the academic level or  in the activities of learned societies alongside their mathematical colleagues. It is quite remarkable how much they could accomplish despite the challenges brought by the war.  Mathematical societies were  established in Moscow, Kiev, Krak\'ow and Lw\'ow  and talks were given. Monographs and textbooks appeared,  journals  were published (``Wiadomo\'sci Matematyczne", ``Prace Matematyczno-Fizyczne", ``Wektor") or planned (``Fundamenta Mathematicae", launch\-ed in 1920). Academic courses went on (with inevitable  interruptions) even though many students and faculty served in the army, buildings were requisitioned  for military purposes and resources (libraries, scientific equipment etc.) were  evacuated.  Doctorates were awarded (Franciszek Leja, Witold Wilkosz, Antoni Plamitzer, Adam Patryn) or recognized (Zygmunt Janiszewski); habilitations were granted (Hugo Steinhaus, Stanis\l aw Ruziewicz, Antoni \L omnicki, Eustachy \.Zyli\'nski, Tadeusz Banachiewicz, Stanis\l aw Le\'sniewski) or denied (Lucjan Emil B\"ottcher). New Polish academic institutions were organized (University of Warsaw, Warsaw Polytechnics, Polish University College in Kiev). Some mathematicians  extended their activity to teaching at a high-school or elementary level when and where the need arose. A few also engaged themselves outside of mathematics and education: in political activities (Wiktor Staniewicz), artistic expression (Leon Chwistek) or writing on cultural and religious Jewish themes (Chaim M\"untz).\\

Did any of  the scholars  mentioned in this article contribute their specific knowledge  to the war effort in the years 1914-18? It does not seem so, at least not directly.\footnote{The cryptography work of Mazurkiewicz, Sierpi\'nski  and Le\'sniewski was done  in the years 1920-21; see \cite{McF} for more information.} Leon Lichtenstein worked for the Siemens company in Berlin, which did support the German military effort, but he dealt with electric cables, not weapons or anything primarily associated with combat. His work, however, was deemed important enough to earn him German citizenship in the early days of the war.\footnote{His fate contrasts with that of Chaim M\"untz, another Russian-born German-educated Polish-Jewish mathematician, who, despite being a serious researcher, educator and  thinker, lost his job at a school in Hessen during the war.}  No other mathematicians mentioned below had ties to industry or institutions of war research. The military service in general did not call for  higher  mathematical skills, although basic knowledge of mathematics and engineering was required sometimes. Edward Stamm served in the Austrian army as a radiotelegraphist, officially translating cablegrams from French, English and Italian, but unofficially might have been involved in deciphering (we have no direct evidence of this, but in 1921 he published a treaty ``On application of logic to the cipher theory"[O zastosowaniu logiki do teorii szyfr\'ow]). Eustachy \.Zyli\'nski had to learn  several engineering subjects in his officer's training in the Russian army in order to become an instructor to future officers. Stanis\l awa Liliental (later Nikodymowa) taught basic mathematics\footnote{We were not able to find a curriculum of these courses. According to \cite{aubingold}, such training was also offered in other armies and included distance measurement and elements of ballistics. Popular texts were \textit{``Soldaten-Mathematik"} by Alexander Wittig (Leipzig 1916) and \textit{``Elementary Mathematics for Field Artillery"} by Lester R. Ford (Louisville, KY, 1919; first circulating as lecture notes).} to Polish army recruits while on leave from her studies at the Warsaw University. While we cannot find  evidence of  anybody's research being directed by war needs,  we can point out several  instances  of mathematical  interests being influenced by war-related circumstances. The most notable cases are those of Stefan Banach, Bronis\l aw Knaster and Kazimierz Kuratowski. All  had to interrupt their studies because of the war --Banach and Kuratowski in engineering, Knaster in medicine-- but they found opportunities for pursuing mathematics and later became pillars of Polish School of Mathematics.  Banach's mathematical career  was spurred by his   serendipitous meeting with Hugo Steinhaus in 1916, which probably would not happen if Steinhaus did not take an administrative job in Krak\'ow after his release from the army. Another interesting case is that of Tadeusz 
Banachiewicz,  an astronomer working in Dorpat, who oriented his research  towards  theory after  the  instruments from his observatory were evacuated deeper inside the Russian Empire. We count him here because the interest in computational methods he developed at that time allowed him later to make some lasting contributions to mathematics. There were other cases of  changing interests, but none of them suggestive of switching from pure mathematics to war-inspired applications.\\

Few individual losses, students developing interests in mathematics as well as arrival  of  a few promising or already established mathematicians (Kazimierz Abramowicz, Jerzy Neyman, Antoni Przeborski, Wiktor Staniewicz, Eustachy \.Zyli\'nski)  after the fall of the Tsarist Russia and the Bolshevik upheaval meant that Poland did not experience a generational gap in mathematics, unlike France.  The so-called Bourbaki thesis (see \cite{aubingold}) claims that the  occurrence of such a gap hindered modern development of mathematics.  Another statement discussed in \cite{aubingold}, Forman's thesis, claims that the war caused collapse of traditional ways of thinking and hence accelerated progress in physics and mathematics, especially in Germany.  We will not argue here whether or how this could be applied to Polish mathematics.  But we have to point out one  rather audacious proposal made during the war which  was crucial for the direction that Polish mathematics took after the war.  It was the publication of the article ``On needs of Polish mathematics"[O potrzebach matematyki polskiej] by Zygmunt Janiszewski in 1918, answering an appeal of a new journal ``Polish Science. Its Needs, Organization and Development. [Nauka Polska. Jej potrzeby, organizacja i rozw\'oj.]", published by  Mianowski Fund. In that article Janiszewski announced his famous program of advancing Polish mathematics by concentrating research on one discipline, possibly in one academic center,  and establishing a specialized scientific journal devoted to this discipline. To a large extent the program was carried out by the Polish School of Mathematics, continuing years after Janiszewski's premature death. However,  foundations were laid before. In 1915, a Polish-language university opened in Warsaw, replacing the Russian one. By conducting lectures and seminars at the University of Warsaw since its very beginning Janiszewski practiced his program before officially formulating it, developing research in topology, aided by like-minded colleagues such as Stefan Mazurkiewicz and  Wac\l aw Sierpi\'nski, and attracting younger  mathematicians to topics of his interest. He understood the importance not only of individual ideas, but also of research collaboration and institutional support. It took many  intellectual and organizational efforts of  Polish mathematicians before the Polish School of Mathematics emerged.  The end of the Great War brought in the independent Polish state, in which the School thrived.

\section{Galicja}

\subsection{Krak\'ow}

Krak\'ow, an ancient Polish capital under Austro-Hun\-ga\-rian occupation, was the seat of the oldest Polish university (established in 1364). In the  academic year 1913/14,   3736  students were attending Jagiellonian University. Classes in mathematics were taught by professors \textbf{Kazimierz \.Zorawski (1866-1953)}\footnote{\.Zorawski, born in the Polish Kingdom, got his doctorate in Leipzig under Sophus Lie. He took a chair of mathematics at the Jagiellonian University in 1893. His research was mainly in differential geometry, mechanics and iteration theory (\cite{zorawski}).}  and \textbf{Stanis\l aw Zaremba (1863-1942)}\footnote{Zaremba, born in the Russian Empire, got his doctorate in Paris. He was appointed a professor of mathematics at the Jagiellonian University in 1900. His research was mainly in differential equations, potential theory and mathematical physics (\cite{dom-zar}).} and docents \textbf{Antoni Hoborski
(1879-1940)}\footnote{In 1911 Hoborski  got habiltation at the Jagiellonian University. Before the war he lectured on descriptive geometry and theoretical arithmetic while also teaching in a gymnasium.  After the war  was nominated for an ordinary professor of mathematics at the newly opened Academy of Mining in Krak\'ow. His main achievements were in differential geometry (\cite{golab69}).}, \textbf{Alfred Rosenblatt (1880-1947)}\footnote{Rosenblatt obtained veniam legendi in 1913. He was unsuccessful in getting a university chair in Poland, so in 1936 he emigrated to Lima, Peru. He worked in the fields of algebraic geometry, analytic functions, mathematical physics and many others. (\cite{CM})} and  \textbf{Jan  Sleszy\'nski (1854-1931)}\footnote{Sleszy\'nski retired from his chair in Odessa in 1909 and moved to Krak\'ow, where he lectured since 1911. He worked on number theory, probability  and logic (proof theory). He used the title of professor, although formally he was appointed to professorship only in 1919. His teaching in the period 1911-19, along with some other classes in mathematics and history of sciences, was financed by the fund of W\l adys\l aw Kretkowski (\cite{dc-kretk}). His chair was the first chair of mathematical logic in Poland, and possibly in the world. \cite{wolenski})}
A course of mathematics for naturalists was taught by \textbf{W\l odzimierz Sto\.zek (1883-1941)}, who also taught in a gymnasium.\footnote{Sto\.zek's doctorate proceedings, started in 1917, were finalized in 1922.}

The first year of the Great War proved very difficult for Krak\'ow, which was a major fortress. The commander Karl Kuk ordered closing the gates on October 17, 1914. The Austrians evacuated over 63 000 people from Krak\'ow. The Russian army advanced fast to the west, coming near Krak\'ow and surrounding it. The first and second battle of Krak\'ow took place respectively in the periods of November 18-22, 1914, and December 6-12, with 40 divisions fighting on both sides.  
Ultimately the Russians were stopped near \L apan\'ow on December 12, 1914, and forced to retreat after heavy losses.(\cite{chwalba})\\

The Jagiellonian University was closed for the winter semester in the academic year  1914/15. Many of the university buildings were requisitioned for military purposes.The Academy of Sciences and Arts, which managed to avoid  requisitions, came to help the University. 
On May 1, 1915,  the University restarted teaching activities after a few month break. It was a success: the classes were held during the shortened summer semester, among them two lectures in mathematics by  Sleszy\'nski.






Regular academic activities continued after 1915, even though some faculty and students were enlisted in the military. \textbf{Marian Smoluchowski (1872-1917)}\footnote{Smoluchowski, an outstanding physicist, made some lasting contributions to mathematics. The most important ones are the theory of Brownian motions and the Smoluchowski equation, which is a limit case of the Fokker-Planck equation.}, who was appointed to the chair of experimental physics at the Jagiellonian University in 1913,  was drafted in 1914 as an Austrian reserve officer. He commanded an artillery detachment, but soon he was released and allowed to come to the University of Vienna, and then back to Krak\'ow. The building of chairs of physics was occupied by the military, so he conducted his work in the former private apartment of his colleague Karol Olszewski (1846-1915).\footnote{Olszewski was the first, jointly with Zygmunt Wr\'oblewski (1845-1888), to liquefy oxygen and nitrogen.} In 1916 Smoluchowski was invited to G\"ottingen to deliver Wolfskehl lectures in June 20-22. His topic was ``Drei Vortr\"age \"uber Diffusion, Brownsche Molekularbewegung und Koagullation von Kolloidteilchen" (Three lectures on diffusion, Brownian motion and coagullation of colloidal particles).  In 1917 he was elected the rector, but  he died of dysentery  before the new academic year began. His duties were taken over by  \.Zorawski 
(\cite{GLM17}, \cite{fulinski98}, \cite{filkrak2}).\\

\textbf{W\l adys\l aw Nikliborc  (1889 - 1948)} passed the wartime maturity exam in  Wadowice on December 14, 1916. In the first semester of the academic year 1916/17  he studied law at the UJ, but soon  (in December 1916) he joined the Polish Legions.  
 He served in the 19th battery of the combat artillery.   After the Oath Crisis he transferred along with the whole regiment to the Austrian army, in which he served until the emergence of the Polish state (he served in France, then in Czech territories in Olomouc). After 2 years he was released from the military and continued his  studies, which he interrupted again twice to take part in the wars with Ukraine and Soviet Russia (1918 and 1920). In 1922 he completed his studies  at  UJ and started working as an assistant to  Antoni \L omnicki at the  Lw\'ow Polytechnics (\cite{nikliborc}).\\


\textbf{Stefan Kaczmarz (1895-1939)}  started studies in mathematics, physics and chemistry in 1913 at the Jagiellonian University. On September 1, 1914 he enlisted in the Polish Legions. He served in the 16th company of the 2nd infantry regiment. From March 1915  to March 1917 he took part in the Carpathian campaign. In July 1917 he was nominated to the rank of  \textit{plutonowy} (platoon sergeant). In 1917 he was transferred to artillery. After the Oath Crisis and dissolution of the Polish Legions he was interned in Huszt  and Bustyahaza (Hungary). He escaped,   but was caught in Galicja and placed in the internment camp in Witkowice (now part of Krak\'ow). In January 1918 he was assigned to the School of Artillery Ensigns of the Polish Auxiliary Corps (Polski Korpus Posi\l kowy) in Walawa (Przemy\'sl region). Released from the army in March 1918, he resumed his studies in Krak\'ow. From November 1918 to February 1919 he served in the Academic Battalion in Krak\'ow as a private. In July 1920 he voluntarily enlisted in the Polish army. He managed to finish his studies in 1922  
(\cite{maligranda07}).\\

Before the war  there were a few Polish paramilitary organizations in Galicja: the Riflemen Union, the ``Rifleman" Society, the Bartosz Troops (Zwi\c azek Strzelecki, Towarzystwo ``Strzelec", Dru\.zyny Bartoszowe). The Austro-Hungarian authorities tolerated them, hoping to use them in a possible conflict with Russia.\footnote{The possibility of such a conflict arose already in 1910, after Austro-Hungarian annexation of Bosnia and Herzegovina.}   Soon they were integrated with the Union for Active Struggle Zwi\c azek Walki Czynnej,   started  by the immigrants from the Polish Kingdom, J\'ozef Pi\l sudski (1867-1935) and Kazimierz Sosnkowski (1885-1969), who fought the imperial Russia in 1905 and intended to  continue the fight. With the outbreak of war, these organizations gave rise to the Polish Legions (Legiony Polskie), a separate formation within the Austro-Hungarian army. In 1917 they were supposed to become a part of Polnische Wehrmacht (Polska Si\l a Zbrojna) organized by the German and Austrian authorities occupying the Polish Kingdom. However, most soldiers and officers refused to swear an oath of loyalty to the German emperor. As a result of the Oath Crisis, the former Russian subjects were interned in the camps in Beniamin\'ow and Szczypiorna, and the Austro-Hungarian subjects were enlisted in the Austro-Hungarian army. The Second Brigade of Legions swore the oath and continue fighting until Austrian command as the Polish Auxiliary Corps (Polski Korpus Posi\l kowy) until the Treaty of Brest-Litovsk  in 1918. In protest  against  territorial concessions to Ukraine, the corps, commanded by General J\'ozef Haller (1873-1960), crossed the front lines to the Russian side. Those unsuccessful were interned in Marmaros-Sziget and Huszt (\cite{kutrzeba}).\\

The  events of war brought to Krak\'ow many people trained in mathematics. \textbf{Edward Stamm (1886-1940)}\footnote{Stamm was a descendant of German Josephine colonists in Galicja.}, a logician, philosopher and historian of science, a graduate of the University of Vienna, taught mathematics in a private gymnasium in Suroch\'ow near Jaros\l aw since 1911. By 1914 he published 23 works on logic and philosophy of mathematics. As the front lines came close, he went first to to Nowy S\c acz, then to Vienna. Many refugees from Eastern Galicja were there, among them Kazimierz Twardowski (1866-1938)\footnote{Twardowski was  a founder of the Lw\'ow-Warsaw philosophical school.} and other Lw\'ow philosophers. Stamm came into contact with Twardowski. Then in 1915 he was drafted into Austrian army and served in a radiotelegraphy station in Krak\'ow, translating cablegrams from French, English and Italian. In 1917 he started officer's training in telegraphy in St. P\"olten and after  completing it,  he commanded  a telegraph station in Cheb (Bohemia). He returned to Krak\'ow in 1918 and served as a Polish commander of the telegraph station until his discharge in 1921, promoted to the rank of a captain (\cite{pabich}, \cite{wachulka}).\\ 

\textbf{Leon Chwistek (1884-1944)}, a mathematician, philosopher  and painter, interrupted his studies of drawing in Paris when the war broke out and joined the First Brigade of Legions. He suffered a wound to his leg, which caused a lasting impairment. After the Oath Crisis Chwistek returned to Krak\'ow, where he devoted himself  to his art and art theory. In 1917 he  co-founded the group of ``Formists" (initially called the ``Polish expressionists"), whose goal was creation of a modern national style, merging the achievements of the Western avant-garde (expressionism, cubism, futurism) with native traditions (medieval arts and crafts, primitive reverse glass painting from Podhale).\footnote{Among the members of the group there were a painter and poet Tytus Czy\.zewski (1880-1945), painters and stage designers brothers Zbigniew Pronaszko (1885-1958) and Andrzej Pronaszko (1888-1961), a sculptor August Zamoyski (1893-1970),  a painter Tymon Niesio\l owski (1882-1965),  a painter, playwright and philosopher  Stanisław Ignacy Witkiewicz (1885-1939)} (\cite{achdaw}, \cite{formisci}, \cite{rzewuski})\\

\textbf{Franciszek Leja (1885-1979)}, after completing his studies in Lw\'ow and obtaining teacher's licence in mathematics and physics,  worked as a teacher in high schools in Krak\'ow and Bochnia. The scholarship from the Academy of Letters enabled him to continue his studies in Paris and London in 1912-13. As a member of Bartosz Troops (Dru\.zyny Bartoszowe)\footnote{Bartosz Troops (Dru\.zyny Bartoszowe) were a military organization active in Lw\'ow since 1908, bringing together the academic youth of peasant origin. In 1914 they were incorporated into the Eastern Legion.}, he was enlisted in the Eastern Legion and fought in 1914-15. When the Legion was dissolved, Leja returned to teaching in the Gymnasium V in Krak\'ow and worked half-time as an assistant in the Chair of Mathematics at the Jagiellonian University at the recommendation of Kazimierz \.Zorawski. Under \.Zorawski's supervision he defended his doctorate in 1916. The title of his thesis was ``W\l asno\'s\'c niezmiennicza r\'owna\'n r\'o\.zniczkowych zwyczajnych ze wzgl\c edu
na przekszta\l cenia styczno\'sciowe [Invariant property of ordinary differential
equations with respect to contiguous transformations.]"\footnote{Leja worked mainly in  potential theory, approximation theory in complex domain and topological groups. He was also a pioneer researcher in several complex variables and is considered a founder of Krak\'ow scientific school of complex analysis.}\\

\textbf{Witold Wilkosz  (1891-1941)}-- Stefan Banach's classmate in the Gymnasium IV in Krak\'ow, showed early interests in mathematics and Oriental languages. A paper on semitology brought him a scholarship to the university of Beirut. Later, in 1912, he enrolled at the Royal University  as an ordinary student at the Faculty of Sciences.  He took courses in mathematics from   Giuseppe Peano (1858-1932), Guido Fubini (1879-1943) and Corrado Segre (1863-1924). Under the direction of   Peano he prepared his PhD thesis concerning the theory of  Lebesgue integral and in 1914 he obtained his doctoral degree. With the eruption of World War I he was called to return to the Austro-Hungarian Empire.\footnote{Italy was a member of the Triple Alliance with Germany and Austro-Hungary, but in May 1915 it revoked the alliance and entered the war on the side of the Allied Powers.} In the first year of the war he fought in the Polish Legions.   Then he continued his studies in mathematics at the Jagiellonian University, which he finished in 1917. It was impossible for him to nostrify his diploma after returning to Krak\'ow, but in 1918 he obtained the degree of the doctor of philosophy for the thesis  ``Z teoryi funkcyi absolutnie ci\c ag\l ych i ca\l ek Lebesgue'a". (On the theory of absolutely continuous functions and the Lebesgue integrals). The supervisor Stanis\l aw Zaremba noted in his report that the thesis made a valuable contribution to modern analysis. Among other things, Wilkosz  corrected a mistake noticed in one of the papers by Charles De la Vall\'ee Poussin (1866-1962). He also included a thank-you note to Banach, reflecting their discussions on mathematics. In the years 1917-18 Wilkosz  taught at private gymnasia in Zawiercie and Cz\c estochowa. He also audited courses in law (\cite{sredniawa-ww}).\\

\textbf{Stefan Banach (1892-1945)} got interested in mathematics as a gymnasium student in Krak\'ow. However, he was not convinced that he could make  a significant contribution to mathematics. So after finishing gymnasium  in 1910 started a course of studies in civil  engineering at the Lw\'ow Polytechnics. In the late spring 1914 he obtained the so-called half-diploma, having passed all compulsory exams for the first and second year of studies as well as the first state licensing exam. After the war erupted in August 1914, the main building of the Polytechnics was requisitioned by the Austro-Hungarian army for a war hospital, and when the Austro-Hungarians retreated in September 1914, it was occupied by the Russian army. 
Banach remained in Krak\'ow until the end of the war, even though the classes at Polytechnics resumed in 1915, after Austro-Hungarians recaptured Lw\'ow. He was rejected in the draft because of being left-handed and having poor vision in his left eye. During the war he first worked as a private tutor in Krak\'ow, then he held a job on road construction, supervising a team of workers.  (\cite{jakim}, \cite{nikon}, \cite{kaluza}, \cite{ciesielscy13})\\ 

It was in the wartime Krak\'ow that Banach embarked on his amazing mathematical career, when Hugo Steinhaus met him by chance one fall day in 1916 (cf. \cite{ciesielscy17}). Here is how Steinhaus recounted his meeting of two men sitting on a bench in the Planty park surrounding the city center and engaged in an advanced mathematical discussion: (\cite{steinhaus})\\
\textit{``Even though Krak\'ow was formally a fortress, one could take walks in Planty in the evenings. During such a  walk I heard the
words `… Lebesgue measure …' -- I went to the bench and introduced myself
to two young adherents of mathematics. They told me that they had another companion, 
 Witold Wilkosz, whom they highly praised. They were Stefan Banach and Otto Nikodym. Since then, we met regularly (...)."} Banach and Steinhaus started collaborating on mathematical problems. After the war Steinhaus arranged for Banach to work as an assistant to Antoni \L omnicki in Lw\'ow. \\

\textbf{Hugo Dyonizy Steinhaus (1887 -1972)}    studied mathematics in Lw\'ow (in 1905 - 1906), then (in 1906 --1911) in G\"ottingen.  In 1911 he  received there PhD degree on the basis of his dissertation entitled Neue Anwendungen des Dirichlet'schen Prinzips. Afterwards he traveled over Europe and published mathematical papers as a private scholar. At the beginning of the war he  moved with his family to Vienna. Then  he went by himself  to Cracow, reported to the  recruiting office of the Polish Legions and 
was assigned to the Military Department of the Principal National Committee (Naczelny Komitet Narodowy, NKN). Initially he did some  office work in  Vienna,  but soon he was sent to to the front, to the 1st Regiment of Artillery of the Polish Legions. He served along with Zygmunt Janiszewski. In 1915 he participated in the war operations in Volyn. His cousin W\l adys\l aw Steinhaus was mortally wounded in a battle and Hugo obtained a leave to attend his  funeral.  He did not come back to his regiment, as his mother managed to  have him recalled from the service. In July 1916 he  took a job in the Center for Reconstruction of The Country (Centrala Odbudowy Kraju) in Krak\'ow. Then he met Stefan Banach, \textbf{Otto
Nikodym (1887-1974)}\footnote{Nikodym taught at the Real Gymnasium IV in Krak\'ow during the Great War.}  and Witold Wilkosz. They started to meet regularly for mathematical discussions in Steinhaus's rented room, joined also by \textbf{W\l adys\l aw \'Slebodzi\'nski (1884-1972)}\footnote{\'Slebodzi\'nski studied at the Jagiellonian University in 1903-1908. In 1913 he went to G\"ottingen, but had to return in 1914 because of the war.  In the years 1919-39 he worked in Pozna\'n. He got his PhD in Warsaw under the supervision of \.Zorawski in 1928. His main research area was differential geometry.}, Leon Chwistek, W\l odzimierz Sto\.zek and Jan Norbert Kr\'oo (1886-?).\footnote{Kr\'oo received PhD in physics in G\"ottingen in 1913. We were not able to find the date of his death.} Steinhaus wrote in \cite{steinhaus} ``...we decided to start a mathematical society", referring to these  meetings. \footnote{This informal society should not be confused with the Mathematical Society in Krak\'ow, which was founded on April 2, 1919.} \\
In 1917 Steinhaus got his habilitation at the Lw\'ow University on the basis of the dissertation ``On certain properties of the Fourier series" (O niekt\'orych w\l asno\'sciach szereg\'ow Fouriera). He arranged to be transferred to Lw\'ow for his job, which he combined with teaching at the university. His lectures on the Lebesgue integral were poorly attended, as nearly all students enlisted in the military. When the Polish-Ukrainian war started in November 1918, he decided to join his parents, in-laws, wife and daughter in Jas\l o. It took him four days to travel the 230-kilometer distance  through Ukrainian and Polish posts as well as no-man's land. Because of his age, he was  exempted from the service in the Polish army. He remained in Jas\l o  working as a mathematical expert at a gas company until  normal activities resumed at the University of Lw\'ow in 1920. (\cite{steinhaus})\\

There was no specialized mathematical society in Krak\'ow until 1919, but mathematicians were active in other scientific organizations, even during the war. For example, on November 29, 1917, Sleszy\'nski gave a talk ``On traditional logic" at the Philosophical Society in Krak\'ow, which he later expanded into a book (published in 1921). Also at the Philosophical Society, on March 1, 1917, Smoluchowski gave a talk ``Remarks on the role of chance in physics" [Uwagi o roli przypadku we fizyce] (\cite{filkrak1}, \cite{filkrak2}.\\

As the Austro-Hungarian empire was collapsing, Krak\'ow  became  the seat of  Polish Commission for Liquidation (Polska Komisja Likwidacyjna), which held  the temporary authority over Galicja and Cieszyn Silesia. It was also the first Polish city to be liberated. On October 31, 1918, a group of Polish soldiers and boy scouts under the command of Lieutnant Antoni Stawarz (1889-1955) took over the railway station in P\l asz\'ow and the army barracks in Podg\'orze and then disarmed an Austrian garrison stationed in the City Hall Tower at the Main Market. The military commander Feldzeugmeister Siegmund von Benigni (1855-1922) handed in power to Polish authorities.  In the book of doctoral promotions at the Jagiellonian University, the words ``Finis Austriae" were entered. Lieutnant Edward Stamm raised Polish white-and-red flag on the radiotelegraph station in D\c ebniki (\cite{jakub}). Civilians aided the military in the effort of keeping the newly gained freedom.  \textbf{Tadeusz Wa\.zewski (1896-1972)}\footnote{Wa\.zewski's achievements were in topology and differential equations. He was a professor at the Jagiellonian University since 1933 and is considered the founder of the Krak\'ow scientific school of differential equations.}, who enrolled as a student at the University in 1915 (first in physics, then in mathematics), served on the Citizen Patrol in Krak\'ow in November and December 1918.\\

\subsection{Lw\'ow} Lw\'ow was the capital of the Kingdom of Galicja and Lodomeria since 1772 (the first partition of Poland).  In 1914 at the  Polytechnic School ,  723 out of 1865 students were studying Civil Engineering, 586  - mechanical engineering, 251 - Technical Chemistry, 243 civil engineering (or architecture), 62 - Engineering Management, which was the first such faculty in Austria. At the University, in the winter semester of the year 1913/14 there were 5871 students enrolled. On September 1, 1914, the city was conquered  by the Russians. Eastern Galicja and Lemkovyna were incorporated into the Russian Empire. Count Georgii Bobrinskii (1868-1928) was made the general-governor of the newly acquired territories.  In the 1914-1915 academic year lectures were not held.  Like in Krak\'ow, some academic buildings were requisitioned for military purposes (e.g., the main building of the Polytechnics served as a hospital). Part of the population was evacuated. Some faculty and students went to Vienna, where teaching and research could be  continued. As a result of breaking the front lines  in May 1915, the Russians withdrew on June 20, 1915, and the Austrians came back on June 23, starting a military rule and repressions against real or perceived supporters of Russia. The academic institutions reopened. The Polytechnic had 130 students in the academic year 1915/16; 198 in   1916/17, 670 in 1917/18, and 989 enrolled for 1918/19. \\

On November 1/2, 1918, Lw\'ow was taken over by the Ukrainians serving in the Austro-Hungarian army, who proclaimed independence of the Western-Ukrainian People's Republic. This started the Polish-Ukrainian war, which lasted until May 22, 1919. Poland reconquered Lw\'ow thanks to the reinforcements from Przemy\'sl under the command of Lieutenant-Colonel (\textit{podpu\l kownik}) Micha\l \ Karaszewicz-Tokarzewski (1863-1964).\\

The chair of mathematics at the university were held by \textbf{J\'ozef Puzyna (1856-1919)} and (until 1918)  \textbf{Wac\l aw Sierpi\'nski (1882-1969)}.\footnote{In 1909 Sierpi\'nski taught the first course in Polish territories on set theory.} \textbf{Zygmunt Janiszewski (1888-1920)} lectured as a private docent. Janiszewski  obtained his doctorate in Paris in 1911 on the basis of the thesis ``Sur les continus irr\'eductibles entre deux points". In 1911-12 he gave lectures in Warsaw at the Society for the Scientific Courses. In 1913 he was nominated to the position of an assistant in the chair of  J\'ozef Puzyna at the Lw\'ow University. Also in 1913, he got his habilitation at the Lw\'ow University on the basis of the thesis `` On dissecting the plane by continua" (O rozcinaniu płaszczyzny przez continua, (published in Prace Matematyczno-Fizyczne 26 (1913), str. 11-63. His habilitation lecture was `` On realism and idealism in mathematics" (O realizmie i idealizmie w rnaternatyce, published in Przegl\c ad Filozoficzny 19 (1916), pp. 161-170).  On  August 30, 1914,  Janiszewski enlisted in the Polish Legions. He took part in the Carpathian campaign in 1914/15. When the Germans took Warsaw, Janiszewski was summoned to a post there. He related later to Hugo Steinhaus that, when he arrived, his driver from the train station to the command was the geometer Max Dehn (1878-1952). \footnote{We were not able to find other accounts confirming Dehn's service in Warsaw; cf. \cite{dehn}, \cite{steinhaus}.}\\
The military service of Janiszewski (then already at the range of a sergeant)\footnote{In 1917, Janiszewski  was the first commander of the men's troop of Polish Military Organization (POW). (\cite{kiepurska}}  and the internment of Wac\l aw Sierpi\'nski in Russia   left Puzyna as the only person  to teach courses (analytic geometry of space and higher analysis) and to run the lower and higher seminar in mathematics in the academic year 1915/16. made it very difficult for the university to offer regular lectures and seminars in mathematics. Because of this difficult situation, in October 1916 the Philosophical Faculty petitioned the Rector to apply for Janiszewski's recall. He was appointed to the post of an assistant for the years 1917 - 1919 and   a request was made to the Imperial and Royal Regency in Lw\'ow  to apply for approval of  this appointment by the Ministry of Religions and Education in Vienna as well as for allocation of funds for his salary. However, Janiszewski did not take this post. After the Oath Crisis he went into hiding to avoid internment. He lived in the Radom region and run (with his own funding) a shelter and a school for homeless children. In 1918 he accepted a chair at the Warsaw University. In the same year he published the article  ``On the needs of mathematics in Poland (O potrzebach matematyki w Polsce),  where he expressed the idea  of creating a mathematical school in Poland. He died in the pandemic of the Spanish influenza.\footnote{It is estimated that 50 to 100 million people died in the pandemic in the years 1918-1920, i.e.,  three to five percent of the world's population. \cite{chwalba}} (\cite{knaster60}, \cite{dom-ksiazka})\\

The nostrification of Janiszewski's doctorate, started in 1914, was delayed because of the war and was finalized only in 1917.  It was necessary for the university to award the right to lecture in order to remedy the staffing shortages. Hugo Steinhaus got his habilitation in 1917. \textbf{Stanis\l aw Ruziewicz (1889-1941)}  obtained  doctorate in philosophy in 1913 at the Lw\'ow University with a thesis ``On a continuous, monotone function which does not have derivative in an uncountable set of points" (O funkcyi ci\c ag\l ej monotonicznej nie posiadaj\c acej pochodnej w nieprzeliczalnej mnogo\'sci punkt\'ow), under the supervision of J\'ozef Puzyna.\footnote{In 1913-14  Ruziewicz received a scholarship from the Academy of Sciences in Krak\'ow from the foundation of W\l adys\l aw Kretkowski. He went to G\"ottingen for the year, which gave him a chance to get acquainted with the problems of contemporary European mathematical research.  In June 1915 he was drafted into the Austrian army. At the beginning he was stationed in Kalusz,  later was commissioned to Hungary and finally to a unit involved in the military censorship in Lw\'ow.  In January 1918, together with the Revision Committee, he was in Bukovina. In February Ruziewicz was assigned to the reserve writers in L\"obnitz, and then was  called again to military censorship in Lw\'ow. He contracted typhus  in L\"obnitz and, in April 1918,  was placed on leave until the end of the year. He stayed in \L a\'ncut, where he taught in a high school.}
In July 1918, while on a leave from the army, he got his habilitation at the Lw\'ow University for the work ``On the monotonic continuous functions with intervals of constancy almost everywhere".  Wac\l aw Sierpi\'nski, back from  his internment, took part in the proceedings (\cite{listy}). Also in 1918,  \textbf{Lucjan Emil B\"ottcher (1872-1937)}\footnote{B\"ottcher got his PhD in 1898 in Leipzig under Sophus Lie (1842-1899). He worked on iteration theory and obtained some pioneering results in holomorphic dynamics.(\cite{dom-st})} made another attempt to get habilitation at the University (he already got one at Polytechnics, where he worked). Like the previous attempts, this one did not succeed. During the wartime (in 1916) one doctorate was awarded at the University: to \textbf{Adam Patryn (1887-1939)}, a gymnasium teacher   in Stryj. The supervisors were J\'ozef Puzyna and Marcin Ernst (1869-1930).\footnote {Ernst was a professor of astronomy at Lw\'ow since 1907.} The title of the thesis was ``Research on functions solving the identity relation $(1-x)^m\Phi(x)+x^m\Psi(x) \equiv 1$" (Badania nad funkcjami rozwi\c azuj\c acymi zwi\c azek identyczno\'sciowy $(1-x)^m\Phi(x)+x^m\Psi(x) \equiv 1$").\cite{prytula}\\
In his evaluation of the thesis, Puzyna wrote:\\

\textit{``Mr. Candidate, in his thesis entitled ``Research on functions solving the identity relation $(1-x)^m\Phi(x)+x^m\Psi(x) \equiv 1$", in his investigations utilized properties in theory of power series, differential equations and combinatorics. He treated the material in a systematic and interesting way, procceding from details to more general cases."}




During the war the following mathematics faculty were active at the Lw\'ow Polytechnics (\cite{dom-ksiazka}):\\

\textbf{Placyd Dziwi\'nski (1851-1936)}: In 1898 he led the 1st Chair of Mathematics at the Polytechnic and held this position until his retirement in 1925. \\

\textbf{Zdzis\l aw Jan Ewangeli Antoni Krygowski (1872 - 1955)}: In 1908 he became an associate professor and in 1909 ordinary professor of mathematics in the Polytechnic School in Lw\'ow, of which he was the rector in 1917 -1918. In the years 1913-1915 he was the dean of the Department of  Water Engineering of the Lw\'ow Polytechnics. In the independent Poland he became a professor at the newly created University of Pozna\'n.\\

\textbf{Kazimierz Bartel (1882-1941)}: After graduating from the Polytechnic School he worked there as an assistant (1907-1911), then Privatdozent of the I Chair of Descriptive Geometry (1911-1912), associate professor (since 1912) and full professor of this chair (since 1917);  in 1912-1939 he led this chair. During the war he served in the Austro-Hungarian army. In 1919 he took part in the Polish-Ukrainian war as the commander of the 1st Railway Battalion. He commanded the defense of the Lw\'ow Main Railway Station. In the independent Poland he held many important political functions, including that of the prime minister.\\

\textbf{Antoni \L omnicki (1881 - 1941)} From 1903 to 1919 he worked as a school teacher in Lw\'ow and Tarn\'ow. In the academic year 1913 - 1914 he was lecturing as a Privatdozent at the Polytechnic in Lw\'ow in the Department of Machine Construction, and  at the Department of Mechanical Engineering at the Lvov Polytechnic.  In 1917 - 1918 he published two works, ``The systems of necessary and sufficient rules for the definition of the concept of quantity"  and  ``On the univalued explicit functions of real variable."  He also got the title of a docent and then became a professor of mathematics, succeeding  Zdzis\l aw Krygowski. In 1918-19 he fought in the Polish-Ukrainian war.\\

\textbf{Antoni Karol Plamitzer (1889-1954)} worked as the assistant at the Polytechnic since 1911. In 1913 he passed the exam for a secondary school teaching licence in mathematics and descriptive geometry. In 1914 he obtained a doctorate degree in technical sciences at the Lvov Polytechnics on the basis of a thesis A Contribution to the theory of flat and curved surfaces under the direction of Kazimierz Bartel (published in ``Wiadomo\'sci  Ma\-te\-ma\-tycz\-ne” 1915, v. 18, 19). \\

On December 3, 1917, Polish Mathematical Society in Lw\'ow had its inaugural meeting.\footnote{Lw\'ow mathematicians were members of other professional organizations before and after World War I. They were particularly active in the  Scientific Society in Lw\'ow (Towarzystwo Naukowe we Lwowie).}   It was established at the initiative of Puzyna, Janiszewski, Steinhaus, \L omnicki, Dziwi\'nski, Krygowski, and \textbf{Tadeusz Cze\.zowski (1889-1991)}, a philosopher whose work had mathematical character.\footnote{Cze\.zowski   taught mathematics and physics in Gymnasium VI in Lw\'ow from 1912 to 1914. In 1915-18 he held a post of the director of the University Chancellor's Office. In 1914 he obtained his doctorate under the supervision of Kazimierz Twardowski  with a dissertation on the Theory of Classes (Teoria klas).   In the independent Poland he first worked for the Ministry for Religious Confessions and Public Education, and in 1923 he took a chair of philosophy at the Stefan Batory University in Wilno (Vilnius).}  The goals of the Society were: support for research in mathematics and related areas, dissemination of mathematical knowledge by scholarly meetings (organized usually every 2 weeks), talks, competitions, publications and collecting the means for learning. The first president was J\'ozef Puzyna. On the board there also were Eustachy \.Zyli\'nski, Antoni \L omnicki and Stanis\l aw Ruziewicz. After Pu\-zy\-na's death in 1919, Marcin Ernst became the president. In the years 1917-1918 he following talks were given at the meetings of the Lw\'ow Mathematical Society: H. Steinhaus: ``Solved and unsolved problems in the theory of Fourier series";
	L. Grabowski ``The harmonic analyzer of Henrici"\footnote{\textbf{Lucjan Grabowski (1871-1941)} studied astronomy and physics in Krak\'ow, Bonn and Munich; a professor  of the Lw\'ow Polytechnics: of surveying in 1909-12 and of spherical astronomy and higher geodesy since 1912.} (in 1917); J. Puzyna: ``On the zero traces of power series”;
	 A. Maksymowicz: ``On Ces\`aro series"\footnote{\textbf{Adam Maksymowicz (1880-1970)}, taught mathematics in Lw\'ow gymnasia and the Polytechnics};
	 Z. Krygowski:  ``On  Tschirnhausen maps in algebra";
	 W. Sierpi\'nski: ``Recent studies on measurable functions";
	 H. Steinhaus: ``On linear and continuous operations in a function field";
	 W. Sierpi\'nski: ``On the continuum hypothesis";
	 W. Sierpi\'nski: ``Definition of the Lebesgue integral without the measure theory";
	 H. Steinhaus: ``Power series in the disk of convergence" (in 1918).

The activities were disrupted by the Polish-Soviet war. In 1920 the Society was dissolved and re-constituted as  the Lw\'ow branch of the Polish Mathematical Society.

\section{The Russian empire}

\subsection{Warsaw}

\par At the beginning of the 20th century two academic-level schools existed in the Polish Kingdom where mathematics was taught\footnote{There were 4 institutions of higher education in the Kingdom in 1914; see \cite{bartn14}.}: the Imperial University in Warsaw  and the  Polytechnic Institute in Warsaw. The language of instruction was Russian.\footnote{Mathematics in the Russian institutions in the Polish Kingdom represented quite a high level.  The most prominent mathermaticians were Dmitri Dmitirevich Mordukhai-Boltovskoi (1876-1952), who worked at the Warsaw Polytechnic, and Georgy Feodosevich Voronoi (1868-1908), who worked at the University. Nikolai Yakovlevich Sonin (1849-1915) spent his entire career at the Imperial University, starting in 1871. Vsevolod Ivanovich Romanovskii (1879-1954), who worked in Warsaw in the years 1911-15, followed as a professor to Rostov-on-Don, and in 1918 to Tashkent. (\cite{duda16}, \cite{nalbaldian})}  A few future Polish mathematicians graduated from the University, e.g. Ka\-zi\-mierz \.Zorawski (in 1888) and Wac\l aw Sierpi\'nski (in 1903).  In the later years some Poles  boycotted these  institutions,  but those who could not go to other provinces of the Empire or abroad  still sought their education  in the Kingdom. \textbf{Zygmunt Chwia\l kowski (1884-1952)} graduated from the Imperial University in 1913 and stayed there to prepare for an academic career. He published a monograph on functional equations in Russian in 1914.(\cite{nalbaldian}, \cite{wbazanekr})\footnote{In the independent Poland Chwia\l kowski taught mathematics at high schools and co-wrote a geometry textbook in 1935 with Wac\l aw Schayer (1905-1959) and Alfred Tarski (1901-1983).  \cite{McF}}\\

At the same time, the Society for Scientific Courses (Towarzystwo Kurs\'ow Naukowych) organized education in Polish at pre-academic and academic level in multiple disciplines, which was not officially recognized, just tolerated. The mathematician \textbf{Samuel Dickstein (1851-1939)}, a graduate of the Imperial University,  was an active promoter of Polish education and scientific organizations. He published at his own expense two mathematical  journals, the first such ones in Polish: ``Wiadomo\'sci Matematyczne"  since 1897 and ``Prace Matematyczno-Fizyczne" since  1888. He  co-founded  the Warsaw Scientific Society (Warszawskie Towarzystwo Naukowe) and donated a library of mathematical books to be used by the Mathematical Study within the Society. In the years 1906-16 he was active in the Mathematical and Physical Circle, which brought together over 100 teachers from the Polish Kingdom.\\
The strife of Poles  for restoration of Polish higher education in the Russian partition and liberalization of education in general\footnote{One of the postulates was admission of women to academic education.} culminated in the massive school strike in the years 1905-08. At that time the authorities made only small concessions, but the break of the war in 1914 brought a mitigation of the Russian policies towards the Polish society. As early as August 14, 1914, the Grand Duke Nikolai Nikolaevich issued an address in which he pledged unification of ``self-governing" Poland under the rule of tsars.\footnote{The address did not have  tsar's authorization.}   With the hopes for freedom  rekindled, the newly established  Warsaw Civic Committee (Komitet Obywatelski  Miasta Warszawy) set up a proposal of restoring the University of Warsaw, \footnote{Dickstein contributed  to this plan.}  which would continue the traditions of the Royal University (1816-1831) and the Main School (1862-1869). Meanwhile, as the Central Powers advanced in the spring and summer of 1915, the Imperial University   was evacuated. The files, libraries and equipment went first to Moscow, where they were followed by the personnel. Then the University moved to Rostov-on-Don, where it remained ever since.\footnote{Even though officially renamed in 1917, it used the name ``Warsaw University" until 1924.} The Polytechnics Institute was also evacuated, to Nizhny Novgorod.\footnote{It later gave rise to the Nizhny Novgorod State Technical University.}\\

On August 5, 1915,  the German army entered Warsaw. The Polish Kingdom was divided into two occupational zones, German and Austro-Hun\-ga\-ri\-an. On September 4, the Germans created the Warsaw general governorate with General Hans von Beseler (1850-1921) as the Governor-General. An idea emerged of creating a Polish  state as one of buffer states  in \textit{Mitteleurope} under political, economical and military  control of Germany.\footnote{The idea was officially announced in a joint declaration by the two respective Governors-General,  von Beseler and Karl Kuk (1853-1935) on November 5, 1916.}  Efforts were made to Polonize the administration and court system. \textbf{W\l adys\l aw (Ladislaus von) Bortkiewicz (1868-1931)}, a professor of statistics at the Friedrich-Wilhelm University in Berlin since 1901 (born  in St. Petersburg to a Polish family and educated there),  was a ``scientific statistical support worker" for the Civil Administration of the General Governorship of Warsaw (Zivilverwaltung des Generalgouwernements Warschau) from November 1916 to February 1917. (\cite{sheynin})\\. 

Opening Polish institutions of higher education was important. Such institutions would prepare future specialists and administrators for the new state in a way that would suit the controlling powers and keep young people out of trouble. Moreover, their existence would improve the attitude of Poles towards the German Empire, so the Central Powers could mobilize Polish men and use resources from the occupied territories.  In these favorable circumstances,  the Civic Committee's project was revisited. The  the Section for Higher Schools (Sekcja Szk\'o\l \ Wy\.zszych) was created, divided into two commissions: the University Commission and the Polytechnic Commission. The mathematician \textbf{Stefan Mazurkiewicz (1888-1945)}, a native of Warsaw,  who studied in Krak\'ow, G\"ottingen and Lw\'ow, was a member of the sub-commission for mathematics and natural sciences. In the fall of 1915 a Polish university and a Polish polytechnic school  were established. Count Bogdan Hutten-Czapski (1851-1937), a Polish aristocrat in German state service, was named a curator, whose function was to act as an official contact  between the General-Government and the administrative structures of the new  schools.  J\'ozef Brudzi\'nski (1874-1917), a physician, became the rector of the University. The rector of the Polytechnic School was \textbf{Zygmunt Straszewicz (1860-1927)},\footnote{Stefan Straszewicz's uncle.} a graduate in mechanical engineering of  Eidgen\"ossische Technische Hochschule in Z\"urich and and a former student of mathematics at the Imperial University,\footnote{Expelled in 1883 in connection with so-called \textit{apuchtinada}.} who until 1916 lectured on differential and integral calculus in the Technical Section of the Society for Scientific Courses and until 1919 taught mathematics and mechanics at the private Mechanical-Technical School of Hipolit Wawelberg (1843-1901) and Stanis\l aw Rotwand (1893-1916). (\cite{bartn15}, \cite{kauffman}, \cite{garlicki}, \cite{duda16}, \cite{kutrzeba}, \cite{chwalba})\\

According to \cite{kauffman},  in the 1915-1916 academic year, the university's teaching staff  included thirty- six lecturers (wyk\l adaj\c acy), the highest rank afforded to teaching staff  at that time, twenty-three assistants (asystenci), and six foreign-language instructors (lektorzy). There were 1,039 students enrolled in 1915-1916. The polytechnic school comprised four departments, where 25 teaching staff instructed about 600 students. \textbf{Kazimierz Kuratowski (1896-1980)}, a Warsaw native who had to interrupt his engineering studies in Glasgow because of the war, was one of the  first students at the University. The introductory courses in mathematics he took   in 1915  were the following:  projective geometry  taught by \textbf{Stefan Kwietniewski (1874-1940)} (in Kuratowski's reminiscence, his lecture was very modern and thorough), analytic geometry by \textbf{Juliusz Rudnicki (1881-1948)}, and  algebra by Samuel Dickstein. 
Along with Kuratowski, the freshman class included, among others,  \textbf{Bronis\l aw Knaster (1893-1980)}\footnote{In 1914 Knaster married Maria Morska (1895?-1945), later a renowned actress, a columnist, and a muse to the Skamander poetic group.\cite{koper}}, who studied medicine in Paris before the war, but in Warsaw switched his interest first to logic (he translated Louis Couturat's \textit{L' Alg\`ebre de la logique} into Polish), then to matematics; as well as \textbf{Stanis\l aw Saks (1897-1942)}, a graduate of the private gymnasium of Micha\l \ Kreczmar in Warsaw.\footnote{The Russians allowed private Polish schools, but their number was limited and they were subject to frequent inspections. After the school strike the situation eased. Kreczmar's school was known for its patriotic atmosphere.} (\cite{kuratow73}, \cite{duda87}, \cite{zyg82}). The lectures in mathematics at the Polytechnic  also included some advanced contents. \textbf{Szolem Mandelbrojt (1899-1983)}, who started his studies there in 1917, recalled Rudnicki presenting Weierstrass' example of a continuous nowhere differentiable function. However, Mandelbrojt found mathematics at the University more attractive, appreciating both lectures and the possibility of private interaction with the faculty, in particular with Zygmunt Janiszewski. He spent two years as a student in Warsaw, published a paper on number theory in 1919, and continued his education in Kharkov, Berlin and Paris. (\cite{mandel})\\

The number of students was rising through years, reaching about 4,500 in 1918. \textbf{Antoni Zygmund (1900-1992)}, who as a gymnasium student  in 1914 was  evacuated with his family  to Poltava, returned to independent Poland  in 1918, completed his education in Kazimierz Kulwie\'c's Gymnasium\footnote{Kazimierz Kulwie\'c (1871-1943), a naturalist and explorer, organized a Polish gymnasium in Moscow in 1915, which he directed for 3 years. After returning to Warsaw he established a school  for  re-immigrants from Russia.} and entered University of Warsaw in 1919.  Zygmund became a student of {\textbf{Aleksander Rajchman (1890-1940)}, a graduate of Sorbonne, who before the war gave private lessons (to Szolem Mandelbrojt, among others), spent the year 1914/15 in Vienna on a scholarship from W\l adys\l aw Kretkowski's\footnote{W\l adys\l aw Kretkowski (1840-1910), a mathematician and a benefactor of science,  graduated from Sorbonne, got a PhD from Jagiellonian University and taught as a private docent at Lw\'ow Polytechnics and Lw\'ow University.} fund (\cite{dc-kretk}), and then worked at the University of Warsaw since 1919, starting at the rank of a junior assistant (\cite{domPBW}).\footnote{Since 1926 Rajchman was a docent at the University and an extraordinary Professor at the Free Polish University (Wolna Wszechnica Polska). He worked on functions of a real variable and probability.} Finally, women gained full access to higher education. \textbf{Stanis\l awa  Nikodym (ne\'e Liliental, 1897- 1988)}\footnote{Otto Nikodym's wife.} enrolled at the University in 1916, but took a break in the year 1918-19 to teach mathematics to army recruits. She returned to her studies and went on to receive her PhD in mathematics in 1925 as the first woman in Poland, to work as an assistant at the Warsaw Polytechnic and to publish several  research papers in mathematics (\cite{dcies-audio}, \cite{stnik}). \\

Many of the faculty had previous  connection to the Society for Scientific Courses. Kwietniewski, who got his PhD at the University of Z\"urich in 1902, concentrated his activities on popularizing mathematics, especially geometry, as well as translating and editing foreign textbooks and monographs. In the years 1907-09 he taught at the Society for Scientific Courses and later contributed to the ``Guide for the self-study". In the independent Poland he continued his university lectures in geometry on a basis of yearly contracts. Rudnicki, a graduate of Sorbonne and later a PhD recipient from Jagiellonian University, taught in private Warsaw schools for men and women since 1912. At the same time he conducted lectures in mathematics at the Society for Scientific Courses (TKN) and at higher pedagogical courses for women. After the Russians' retreat from Warsaw he was a member of the Polytechnic Commission of the TKN  (the electrotechnical-mechanic sub-commission). After creation of the Warsaw Polytechnic he taught mathematics and (briefly) physics there. He was also active in the Society for Aid to the Victims of War (Towarzystwo Pomocy Ofiarom Wojny). Later he became a professor in Vilnius University. (\cite{krolik89})

Among the faculty of the University there also were Stefan Ma\-zur\-kie\-wicz  and \textbf{Zygmunt Janiszewski (1888-1920)} (later the recipients of the first two chairs in mathematics). This was remarkable on two accounts. First,  both of them were rather young; second, both  were  previously connected to Lw\'ow (Janiszewski taught there as a docent and Mazurkiewicz got his doctorate under Sierpi\'nski in 1913)  and Beseler (contrary to the wishes and efforts of the Civic Committee) did not want too many professors from Austrian universities coming to Warsaw.\footnote{According to Kuratowski,  \cite{kuratow79}, Mazurkiewicz  supported the candidacy of the Polonized archduke Karol  Stefan Habsburg-Lotary\'nski (1860-1933) for the throne of the Polish Kingdom under the auspices of Austro-Hungary and Germany.} However, they had significant mathematical achievements and scholarly output:  Mazurkiewicz published 16  works and Janiszewski 20 works before 1916.  Janiszewski was available only part-time, as he was still serving in the military. Since 1916 Janiszewski  and Mazurkiewicz conducted a seminar in topology at the university, possibly the first one in the world in this new discipline. Kuratowski, Knaster  and Saks participated as  students. (\cite{kuratow79}). \\

Student organizations played a major role in the academic life. For example, the analysis of Stanis\l aw Zaremba's ``Theoretical Arithmetic" by \textbf{Jan \L ukasiewicz (1878-1956)}\footnote{\L ukasiewicz was a philosopher and a logician whose work was taking more and more mathematical character. He was a pioneer of multi-valued logic and an inventor of the Polish notation.}  in his course on methodology of deductive sciences prompted a discussion involving professors and students. The discussion of related issues continued  in the meetings of Mathematical and Physical Circle of the Warsaw University, with Kuratowski giving a two-part talk on December 6 and 13, 1917, ``On the definition of a quantity", which soon became his first scholarly publication (see \cite{kuratow79}). On a larger scale, events related to Polish history and Polish national heroes were commemorated. However, the arrest of two students after celebrating the anniversary of the 3rd of May Constitution, led to the students' strike in 1917 and temporary closing of the schools. The occupational authorities transferred the control of the schools to the Temporary Council of State, the first government of the Kingdom of Poland. J\'ozef Pi\l sudski, who held the authority over the military mattters, resigned from the Council, which led to the so-called Oath Crisis. The Council was ultimately disbanded in August 1917. In the fall of 1918 Polish army started to organize itself  and students were joining in great numbers. The Academic Legion (Legia Akademicka)-- a unit comprised entirely of students was formed.\\

The creation of Polish academic institutions did not eliminate the need for the Society for Scientific Courses. There still were many people aspiring to higher education with insufficient credentials for admission, so the Society  
continued its activity during the war. The Department of Mathematics and Physics separated from the Department of Mathematics and Natural Sciences in 1915.  In 1915/16  it run the following compulsory courses: descriptive and projective geometry, taught by \textbf{Wac\l aw Gniazdowski (1864-1938)}\footnote{Gniazdowski, a textile engineer,   taught mathematics and technological subjects at the Technical School of the Warsaw-Vienna Railway (Szko\l a Techniczna Drogi \.Zelaznej Warszawsko-Wiede\'nskiej).  After the school's evacuation he founded his own private 7-grade school of technology and transportation (\cite{kolejowa}). He also taught mathematics at the Real 7-Grade School directed by Witold Wr\'oblewski (1839-1927) with instruction in Polish in the years 1915-18. Later he was a docent at the Warsaw Polytechnics, teaching principles of perspective at the Department of Architecture. (\cite{gniazd})}; analytic geometry, taught in the fall by \textbf{Romuald Witwi\'nski (1840-1937)}\footnote{Witwi\'nski authored several papers and problem books in geometry.} and in the spring by \textbf{Tadeusz Gutkowski (1881-1962)}\footnote{Tadeusz Gutkowski, an optical engineer, graduate of Institut d'Optique in Paris, taught mathematics in Warsaw high schools, and later worked in the optical industry. (\cite{gut}}; introduction to analysis by \textbf{W\l adys\l aw W\'ojtowicz (1874-1942)}\footnote{W\l adys\l aw W\'ojtowicz -- an editor of the journal ``Wektor" for teachers, editor of a series published by the Mianowski Fund, author of high school geometry textbooks and logarithmic tables for the school use.\cite{dom-zar}}. The optional courses in the first semester were  differential and integral calculus, taught by Juliusz Rudnicki, and vector  calculus, by \textbf{Wac\l aw Werner  (1879-1948)}.\footnote{Wac\l aw Werner  studied electrotechnology in Darmstadt, mathematics and physics in Krak\'ow, G\"ottingen and Fribourg. He received a doctorate from the Faculty of Mathematics and Natural Sciences in Fribourg. In 1909-39 he taught physics in  high school  in Warsaw. In 1916/17 he was the dean of the Department of Mathematics and Physics in the Society for Scientific Courses. During that time he co-managed family-owned metal works. Since 1915 he worked at Warsaw Polytechnic, lecturing and conducting experiments; named a titular professor in 1948.\cite{werner}} Lectures by Stanis\l aw Le\'sniewski and Stefan Mazurkiewicz were also planned, but did not run. In the first and second semester there were respectively 27 and 19 students. In 1916/17, \textbf{Franciszek W\l odarski (1889-1944)},  a geometer with doctorate from the University of Fribourg, started to lecture. \\
Mathematical subjects were also taught at the Department of Technology, among them trigonometry by \textbf{Tomasz \'Swi\c etochowski}\footnote{\'Swi\c etochowski taught mathematics at the Real 7-Grade School directed by Witold Wr\'oblewski with instruction in Polish in the years 1915-1919.\cite{czacki}},  algebra with geometry by \textbf{Bruno Winawer (1883-1944)}\footnote{Bruno Winawer-- a physicist, writer and popularizer of science, a graduate of University of Heidelberg.} and analytic geometry by \textbf{Lucjan Zarzecki (1873-1925)}\footnote{Lucjan Zarzecki--a mathematician and educator, a graduate of St. Petersburg University in 1897.}. The recitation classes were taught by \textbf{F. \L a\-zar\-ski}\footnote{We were not able to find the full first name and the dates of birth and death for him.} (differential and integral calculus), R. \'Swi\c etochowski\footnote{Born 1882 (\cite{kiepurska}); probably a misprint of the initial.} (descriptive geometry), A. Winawer\footnote{Probably a misprint of the initial.} (high school mathematics), W. W\'ojtowicz (higher mathematics as well as analytic geometry, together with the lecture). 
 Later the Society also gave rise to the Free Polish University (Wolna Wszechnica Polska), a fully accredited private university operating in the years 1918-1952 in Warsaw and \L \'od\'z 
(\cite{maligranda17}, \cite{10lecie}).\\

The publication of the series ``Guide for the self-study" (Poradnik dla samouk\'ow)\footnote{The series appeared in several cycles in the years 1898-1932, financed by the Mia\-now\-ski Fund. It was meant as an educational aid at an academic level. Each  volume presented the development and the state-of-the art of a given scientific topic, along with exhaustive bibiography. The editors were Aleksander Heflich (1866-1936) and Stanis\l aw Michalski (1865-1949).} continued during the war. A volume on mathematics, starting the second series, was published in 1915. It contained chapters written by Jan \L ukasiewicz (On Science), Zygmunt Janiszewski (General Introduction; Introduction to Level III; Ordinary  Differential Equations;  Functional, Difference and Integral Equations; Series Expansions; Topology; Foundations of Geometry; Logistics; Philosophical Issues of Mathematics; Conclusion; Informational Section), Stefan Kwietniewski (Level I; Level II; Methodology of Teaching; Analytic Geometry; Synthetic and Descriptive Geometry; Differential geometry; History of Mathematics: History of Mathematics in General; History of Mathematics in Poland), Wac\l aw Sierpi\'nski (Arithmetics; Number Theory; Higher Algebra; Set Theory; Theory of Functions of a Real Variable; Differential and Integral Calculus; Differential Calculus and  Summation), Stanis\l aw Zaremba (Theory of Analytic Functions; Differential Equations with Partial Derivatives; Theory of Groups of Transformations; Calculus of Variations) and Stefan Mazurkiewicz (Theory of Probability) (\cite{zpb-poradnik}). Marian Smoluchowski contributed a chapter on physics to the 1917 volume.\\

\subsection{St. Petersburg (Petrograd)} Polish presence was very strong in the Russian capital. In 1910 the number of Poles living there reached its historical maximum of about 65 000 (3,4\% of the total population of the city) . Polish nationals could be found among officers, civil servants, artists and scholars. (\cite{garczyk}) \textbf{Julian Karol Sochocki (Yulian Vasilievich Sokhotsky, 1842-1927)}, born in Warsaw, was educated at the University of Saint Petersburg and was a professor of mathematics there. His results in the field of one complex variable (Sochocki-Casorati Weierstrass theorem, Sochocki-Plemelj formula) became classic. \textbf{Jan Ptaszycki (1844-1912)} was a professor of mathematics at the University of Saint Petersburg and at the M\-ikhai\-lov\-ska\-ya Military Artillery Academy.  His work dealt with elliptic functions and algebraic differentials. \textbf{Wiktor Emeryk Jan Staniewicz (1866-1932)}, born in Samara, educated in St. Petersburg, held the chair of mathematics at St. Petersburg  Polytechnic Institute since 1902. He worked in number theory and mathematical analysis. In 1909 his state service was   suspended for three years, because illegal political activities were discovered to go on in dormitories that he supervised. In that period he lectured as a contract professor.  In the years 1915-17 he was the dean of the Faculty of Civil Engineering, in 1917-18 a vice-rector. Polish mathematicians in St. Petersburg did not form a separate learned society, but were active in the Polish Union of Physicians and Naturalists (Zwi\c azek Polski Lekarzy i Przyrodnik\'ow). Sochocki also presided over St. Petersburg Mathematical Society in the years 1984-1927. (\cite{domPB})\\

The year 1917 brought dramatic political and social changes to the Russian Empire. The (changing) authorities were trying to win the support of Poles. In December 1916  Tsar Nicholas II as the commander in chief issued an order number 870 to land and maritime armed forces, which among the goals of further campaign mentioned the ``creation of free Poland". In March 1917 the Provisional Government stated that it counted on forming a ``free military union" with Poland in the future, while the Petrograd Soviet of Workers' and Soldiers' Deputies (later taken over by the Bolsheviks) announced Poland's right to complete political independence, and the general right of nations to ``self-determination" (\cite{Za17}). In the circumstances favorable to the Polish cause, in July 1917 Staniewicz became a president of the Polish Radical-Democratic Union in Lithuania and Belarus (Polski Zwi\c azek Radykalno-Demokratyczny na Litwie i Bia\l orusi) and took part in the attempts to form the Polish National Executive Commission (Polska Narodowa Komisja Wykonawcza) in Russia. In October 1919 he moved to independent Poland and became a professor of mathematics at Stefan Batory University in Wilno (Vilnius). He was the first president of the Polish Mathematical Society (\cite{wstan}, \cite{iw}).\\

\subsection{Moscow}

According to the census from 1897, there were 9236 Poles living in Moscow, 0.89\% of its population. The massive evacuations from the Polish Kingdom at the beginning of war raised this number.\\

 \textbf{Boles\l aw M\l odziejowski (Boleslav Kornel'jevich Mlo\-dze\-ev\-skii, 1858-1923)} was born in Moscow in a physician's family.\footnote{M\l odziejowski's maternal grandfather was a Czech musician Vincenz (Vikientii) Lemoch (1792-1862?), a brother of Wojciech Ignacy Lemoch (1802-1875), who was a professor of geometry and a rector of the Lw\'ow University. (\cite{ilemoch}, \cite{vlemoch})} He graduated from Moscow University in 1880 and became a professor of mathematics there in 1892. His research interest was in geometry.  In 1902 he served as an opponent in the doctoral defense of Antoni Przeborski.\footnote{Other faculty refused, citing anti-Polish regulations from 1864.\cite{odyniec}}. In 1911 he resigned from his position in protest against decisions of the enlightment minister Lev Aristidovich Kasso (1865-1914).\footnote{Kasso proposed  new ways of staffing vacant chairs at Russian academic institutions, which met with disagreement of the professors.} He continued his teaching activities at the Higher Courses for Women as well as at  Moscow City People's University.\footnote{The People's University was also referred to as Shanyavskii's University, after its founder Alfons  Lvovich Shanyavskii (Alfons Fortunat Szaniawski; 1837-1905), a general of Polish origin in the Imperial army.  It was a research-oriented university, open to anyone regardless of their origin, education, gender, age, nationality, or religious beliefs.  It operated in the years 1908-1918. (\cite{rag})}   At the latter, he conducted lectures in differential geometry and introduced modern-style seminars. In January 1914, he chaired the organizing committee of the Second All-Russian congress of lecturers in mathematics, in which 20 speakers from Polish territories took part.  After the February revolution in 1917 M\l odziejowski  returned to the university.  In 1921 he became the first director of the newly created Research Institute in Mathematics and Mechanics at the Moscow University. (\cite{ZvPu}, \cite{odyniec})  \\

\textbf{Stanis\l aw Le\'sniewski (1886-1939)}, born in Serpukhov in the Moscow governorate and brought up in Irkutsk, studied philosophy and mathematics in Germany, Switzerland and Russia. He completed his doctorate in philosophy in 1912 at the Lw\'ow University under the direction of Kazimierz Twardowski. Afterwards he taught at a school in Warsaw. When the war broke out, he went to Moscow.\footnote{It is not clear why Le\'sniewski went there.  It could be as a result of an evacuation or in connection with his political activities as a member of Social Democracy of the Kingdom of Poland and Lithuania (Socjaldemokracja Kr\'olestwa Polskiego i Litwy).}  He taught matematics at a Polish gymnasium and at the Real School of the Polish Committee of Aid to the War Victims (Szko\l a Realna Polskiego Komitetu Pomocy Ofiarom Wojny), founded for boys from families that were evacuated from the Polish Kingdom.\footnote{The school was directed by the distinguished educator W\l adys\l aw  Gi\.zycki (1875-1947). Among the students there was the future poet Konstanty Ildefons Ga\l czy\'nski (1905-1953) (\cite{galczynska}).} Le\'sniewski was also active in the Polish Scientific Circle (Polskie Ko\l o Naukowe). Through the Circle, he published his book ``Foundations of general set theory, part I" in 1916. The second part was planned for 1917, but  never came out. Despite the title, the book treated Le\'sniewski's own  theory of parts, wholes and concrete collections, which was later  developed into his system of Mereology. In 1918 he returned to Warsaw and on December 14 he submitted  his habilitation dissertation in logic and philosophy of mathematics to be evaluated by Wac\l aw Sierpi\'nski. \footnote{The dissertation consisted of the works ``Problems of the General Theory of Sets, I" and ``A Criticism of the Logical Principle of the Excluded Middle".}   In 1919 Le\'sniewski became 
a professor of philosophy of mathematics at the University of Warsaw. In 1920, along with Stefan Mazurkiewicz and Wac\l aw Sierpi\'nski, he contributed to breaking Soviet codes in the Polish-Soviet war. (\cite{SLstanford}, \cite{McF},  \cite{betti})\\

\textbf{Wac\l aw Sierpi\'nski (1882-1969)}-- a native of Warsaw, a graduate of the Imperial University (under the direction of Georgy Voronoi) and a PhD recipient from the Jagiellonian University, was an extraordinary professor of mathematics at the Lw\'ow University since 1910. The outbreak of the war found him in Belarus, in the estate of his parents-in-law. As an Austro-Hungarian citizen, and hence an enemy alien, he was interned in the city of Vyatka (nowadays Kirov). Thanks to the efforts of Moscow mathematicians (mainly Dmitri Fyodorovich Egorov, 1869-1931) he was allowed to relocate to Moscow in 1915. The Rectorate of the Lw\'ow University was notified of Sierpi\'nski's internment in Moscow through the American consulate\footnote{The United States of America remained neutral in the war until April 6, 1917.} in Vienna in February 1916. The university administration made efforts to transfer to Sierpi\'nski his overdue (since 1914) salary using the same diplomatic channels, but they were unsuccessful.\\

The Moscow period was very fruitful for Sierpi\'nski.  It marked a beginning of his deep studies of the axiom of choice and its role in mathematics. He gave a talk on the subject at the meeting of the Moscow Mathematical Society on February 21, 1917.\footnote{It was preceded by a note in  ``Comptes Rendus" of the French Academy in 1916. The expanded version of the talk was later published in French as ``L'axiome de M. Zermelo et son r\^ole dans la th\'eorie des ensembles et l'analyse" in ``Bulletin International de l'Acad\'emie des Sciences de Cracovie. Classe des Sciences Math\'ematiques et Naturelles, S\'erie A" 1918, s. 97-152 and in Russian as as ``Aksioma Zermelo i eio rol' v teorii mnozhestv i analize", Matematicheskii Sbornik 1922, tom 31:1, 94-128. See also \cite{lewand}.}   He also started a collaboration and friendship with Nikolai Nikolaevich Lusin (1883-1950).\footnote{The scientific relations later also extended to Lusin's students, who visited Poland and published in Polsih journals.} In the years 1915-1918 he published 41  papers, 4 of them jointly with Lusin and 3 on problems related to Lusin's research (\cite{sierp}).
While interned, Sierpi\'nski was active in the Polish Scientific Circle (Polskie Ko\l o Naukowe) established in November 1915 in Moscow. Through the Circle  he published the first volume of ``Mathematical Analysis", which 
he dedicated to the Polish University in Warsaw.\footnote{Like in the case of Le\'sniewski, the second volume was planned, but never appeared.} He also gave talks in the Moscow Mathematical Society. In February 1918 he returned to Poland through Finland and Sweden. He resumed his lectures in Lw\'ow in the summer semester  1918, but in the fall he moved to Warsaw. He was nominated for an ordinary professor of mathematics at the philosophical faculty of the Warsaw University  by the decree of the Chief of State from March 28,  1919. He announced his resignation from the chair in Lw\'ow in a letter  dated May 19, 1919, thanking his colleagues for the kindness they offered him during his stay in Lw\'ow. (\cite{mioduszewski}, \cite{sinkiewicz}, \cite{sierp-teka}).\\




\textbf{Kazimierz Jantzen  (1885-1940)}  got a doctorate in astronomy in Munich in 1912. In the years 1912-14 he was at the astronomical observatory in Potsdam. The outbreak of the war found him in Warsaw. As a German citizen, he was  interned by the Russian authorities in Vyatka, and then transferred to Moscow. He taught in Polish high schools and was active in the Polish Scientific Circle. He published a book ``On the influence of the spectral type of stars on determining the apex of the Sun" (O wp\l ywie typu widmowego gwiazd na wyznaczanie apeksu s\l o\'nca). He returned to Poland in 1918, worked at the astronomical observatory in Warsaw, Warsaw Polytechnics (lecturing on advanced surveying and the error theory) and the Military Geographical Institute. Then he took a chair of astronomy at the Univeristy of Wilno (Vilnius), where he also lectured on analytic geometry, statistics amd mathematics for naturalists. (\cite{dom-koron}, \cite{rybka})\\

\subsection{Kharkov}

The University of Kharkov was established in 1805 by a Polish aristocrat, Seweryn Potocki (1762-1829). In 1897 there were 3969 Poles living in Kharkov,  2,28\% of its population. 
  The distance  from the front lines of the World War I allowed for the university activities to go on as usual, at least at the beginning of the war. \textbf{Jerzy (Yuri Cheslavovich) Neyman (1894-1981)}, born in Bendery, in the Bessarabia governorate, entered his second year of studies at the University of Kharkov in 1914/15.  Rejected in the draft because of poor eyesight, he was preparing a paper on Lebesgue integral to enter a university-wide competition at the encouragement of his professor \textbf{Cezary Russjan (1867-1935)}.\footnote{Russjan got his doctorate in Odessa in 1900. For some time he held lectures at the Lw\'ow Polytechnics and the Jagiellonian University. In 1907 he took the chair of mathematical analysis at the University of Kharkov. His main interests were differential equations and probability.} His 580-page work won and he received a monetary equivalent of the gold medal (the actual medal could not be awarded because of wartime restrictions on  metals). In 1917 Neyman finished his course of studies and, on Russjan's recommendation, was granted a government stipend to prepare himself for an academic career. At the same time he started working in the Kharkov Institute of Technology as an assistant to \textbf{Antoni Przeborski (1871-1941)}\footnote{Przeborski got his doctorate in 1902 at Moscow University. In 1908 he became an ordinary professor at the University of Kharkov. He also taught at Kharkov Polytechnic Institute, Women's Higher Courses and the Workers' University. His main interests were analytic functions, dfferential equations and variational calculus. In the independent Poland he was a professor of universities in Warsaw and Vilnius; he also taught mathematics and mechanics at the Warsaw Polytechnics.}, in analytic geometry and introduction to analysis,  and as a lecturer in elementary mathematics. \\

The wartime situation in the Russian empire was  complicated by the outbreak of 2 revolutions: in February and October 1917, and by the Ukrainian-Soviet war.  Ukrainian People's Republic of Soviets  formed in 1917 in Kharkov  fought Ukrainian People's (or National) Republic  proclaimed in January 1918 and based in Kiev, which was aided by the Germans after  the Brest-Litovsk peace treaty between Germany and Russia was signed on February 9, 1918. The fightings continued after the Germans withdrew. The University and Polytechnic Institue in Kharkov continued to operate under the Bolshevik rule (with some interruptions), opening their doors to many more people from underprivileged background. Neyman was assigned a task of teaching remedial classes in mathematics to these new students. In addition, he taught mathematics and physics in a newly opened Polish high school in Kharkov. He also spent a brief  time in prison, arrested for bartering matches for food in the black market. 
It was during the wartime that Neyman   got interested in statistics (which later became the field of his highest achievements),  through discussions he had with Sergei Natanovich Bernstein (1880-1986),  newly promoted to professorship at the University of Kharkov. \\

In the years 1919-1920 Przeborski was the rector of the university, reorganized into the Academy of Theoretical Sciences.  Neyman recalled that  in the time of severe deprivations Przeborski arranged for the professors to obtain permission to chop trees in the nearby park for fuel. Due to a misunderstanding about  legitimacy of the permission several professors were arrested, including Przeborski himself. He was released, later made the dean, and then the rector again. Eventually, Neyman and Przeborski left for Poland after the Polish-Soviet war (in 1921 and 1922, respectively).  Russjan remained in Kharkov until his death in 1935.\footnote{He was dismissed from his university position  by the Soviet authorities in 1934.}(\cite{reid}, \cite{kijas})\\

\subsection{Kiev} The Imperial St. Vladimir Kiev University, established in 1843, could be viewed as contination of  the Krzemieniec Lyceum, since it started with the Lyceum's  assets and its Polish faculty. It was a popular destination for Polish students; the poets Boles\l aw Le\'smian (1877-1937) and Jaros\l aw Iwaszkiewicz (1894-1980) studied there.  Overall,  there were 16 579 Poles in Kiev in 1897, 6,69\% of its population. In 1919 the number went up to 36 800, or 6,77\%. \textbf{Eustachy \.Zyli\'nski (1889-1954)},  born to a Polish family in the Podolia district, graduated from the St. Vladimir University  as a student of Dmitrii Alexandrovich Grave (1863-1939). He passed his exams for master's degree in 1914 in Kiev after taking a study trip to G\"ottingen, Cambridge and Marburg, then completing his exams and presenting a thesis ``On the field of $p$-adic numbers". From 1912 to 1915 he worked at the St. Vladimir University. On April 16, 1916 he was drafted into the Russian army as a \textit{Praporshchik} (ensign). As part of his officer's training, he completed several courses in engineering subjects and in radiotelegraphy in Kiev and Petersburg.  On February 7, 1917, he was nominated to the rank of \textit{Podporuchik} (Second Lieutenant). He became the Head of Radiotelegraphy of the South-Western Front, then he commanded an officers' class. He did not engage in combat; his main task was to teach electrotechnical subjects to the Staff of 103rd Front in Kamenets-Podolskii (Kamieniec Podolski; Kamianets-Podilskiy) and Berdichev. \\

On July 24, 1917, Polish I Corps was formed in Belarus from Polish soldiers serving on Northern and Western Fronts, under the command of General J\'ozef Dowb\'or-Mu\'snicki (1867-1937). \.Zyli\'nski reported to the commander of the Corps in November 1917.  In the period from December  1917 to February 1919 he worked  in the Polish University College, Ukrainian State University and Higher Private Polytechnic Institute in Kiev. He taught classes in analytic geometry, set theory,   higher algebra and introduction to analysis. He got his habilitation  at the Polish University College. He published one paper, ``O zasadach logiki i matematyki" (``On the principles of logic and mathematics")  in the Reports of Polish Scientific Society  in Kiev (Sprawozdania Polskiego Towarzystwa Naukowego w Kijowie) in 1918. He also wrote 2 extensive works in the fields of  algebra and logic during the period of war (both probably remained unpublished).  On February 19, 1919, \.Zylinski went to Warsaw as an officer in the Polish army. He taught in the Officers' School of Communication. He remained in the military service until September 1919, still as a second lieutenant  (\textit{podporucznik}) in a radio-telegraph battalion. He was released to 
become an extraordinary professor of mathematics at Jan Kazimierz University in Lw\'ow.\footnote{\.Zyli\'nski's candidacy was supported by a mathematical committee, consisting of an astronomer Marcin Ernst (1869-1930), a physical chemist Roman Negrusz (1874-1926), a philosopher KazimierzTwardowski (1866-1938) and a physicist Ignacy Zakrzewski (1860-1932).} Earlier, he rejected an offer to take a chair of mathematics 
 at the Kamianets-Podilsky State Ukrainian University (formed in 1918 under a law signed by Pavlo Skoropadskyi  (1873-1945), Hetman of Ukraine).  In  Lw\'ow he soon became the head of the Chair A. He initiated a revival of algebra in Lw\'ow (\cite{maligranda09},\cite{domSZ}).  \\

\textbf{Kazimierz Abramowicz  (1888-1936)}, born in the Polish Kingdom, finished his course of studies in mathematics at St. Vladimir's University in Kiev in 1911, receiving a gold medal for his work ``On hypergeometric functions with one removable singular point". He worked under direction of Boris Yakovlevich Bukreev (1859-1962). In 1914 he passed his master's degree exams and went to Berlin and G\"ottingen for further studies. Because of the outbreak of the war, he returned to Kiev. In the academic years 1914/15 and 1915/16 he lectured at the Kiev Polytechnic Institute. He was delegated as a \textit{docent} to the Perm branch of the Petrograd University for the year 1916/17.\footnote{The  branch   was established in 1916 as a result of evacuation of the Petrograd University  deep into the territory of the Empire, to safeguard people and to alleviate provisional shortages. Mathematicians Yakov Davidovich  Tamarkin (1888-1945), Alexander Alexandrovich  Friedmann (1888-1925), Abram Samoilovich Besicovitch (1891-1970), Nikolai Maximovich Gjunter (1871-1941), Rodion  Osievich Kuzmin (1891-1949) and Ivan Matveevich Vinogradov (1891-1983)  taught  there in the early years.  (\cite{demidov15})} As the branch  became an independent university in 1917, Abramowicz was nominated an extraordinary profesor in the chair of mathematics. Because of the war operations  he could not return to Perm in the fall of 1918, so he taught recitation classes in mathematics at the Polytechnic Institute in Kiev. He returned to Poland in June 1920 and started working at the newly established University of Pozna\'n  in 1921. (\cite{maligranda16}) \\

\textbf{Izabela Abramowicz (1889-1973)}\footnote{Kazimierz Abramowicz's sister.} was the first woman to receive the 1st degree diploma at the Faculty of Mathematics and Physics of the St. Vladimir University in Kiev and a gold medal for the thesis ``On double integrals on algebraic surfaces.". Like her brother, she worked under direction of Boris Bukreev. She stayed at the university, by permission of the education minister, but without a stipend,  to prepare herself for exams towards her master's degree. She also taught at three gymnasia in Kiev. 
In the years 1917-1920 she lectured on introduction to mathematics (as a \textit{docent}) at the Polish University College in Kiev.\footnote{The College was established at the initiative of Wac\l awa Peretiatkowiczowa (1855-1939), a headmistress of two women gymnasia in Kiev. It started in 1916 as Higher Polish Learning Courses, initially allowed  to offer only a program in humanities. It  continued in 1917-1919 as the Polish University College. The faculty was recruited from Poles teaching at Russian institutions of higher education as well as academics from the Polish Kingdom and Galicia who for various reasons found themselves in the Russian Empire. About 40 people overall taught there. Some of them, including Eustachy \.Zyli\'nski, obtained their habilitation at the College. The  students-- mainly Polish nationals from the Kiev Governorate, along with some incomers from the Polish Kingdom or Galicia-- were interested in getting higher education and preparing themselves for professional specialization or teaching in Polish schools in the Western Ukrainian territories. The majority  of them-- 572 out of 718 in the first semester-- were women. The College's activities were financed mainly by the Society for Supporting Polish Culture and Learning in Ruthenia (Towarzystwo Popierania Polskiej Kultury i Nauki na Rusi). (\cite{RZ})} She joined the College when it expanded its course offer  to mathematics and sciences.    She was one of two women among the  faculty members; the other one was  Antonina Dylewska (1883-1951), a mineralogist. In addition to her teaching activities, Abramowicz was also a member of the short-lived Polish Scientific Society in Kiev. The College faculty and students started to leave as  the fightings  continued. Even after the Great War had ended, Kiev changed hands, passing from Germans to Ukrainians to Bolsheviks to White Russians to Ukrainians and Poles to Bolsheviks again.  Abramowicz was supposed to leave Kiev in 1920 with the retreating Polish army. In August 1923 she arrived to Pozna\'n. During the Second Republic and after World War II she taught mathematics in high schools in Pozna\'n. (\cite{maligranda16})\\

\subsection{Yuryev (Dorpat; Tartu)} The University of Dorpat continued the traditions of a Jesuit college established by the Polish king Stefan Batory in 1583 and attracted many Polish students.\footnote{Until 1893 the language of instruction  was German. Then the university was fully Russified. The city itself was renamed Yuryev.} Among distinguished graduates there were Tytus Cha\l ubi\'nski (1820-1889), a physician and promoter of tourism in Tatra mountains, and Wincenty Lutos\l awski (1863-1954), a philosopher. \textbf{Tadeusz Banachiewicz (1882-1954)}, a native of Warsaw,\footnote{Banachiewicz studied at Warsaw, Kazan, Moscow and G\"ottingen.  He was involved in the activities of the Society for Scientific Courses. He was primarily an astronomer, but made a lasting contribution to mathematics by inventing the (non-associative) algebra of Cracovians. He is also credited with a proof of Schur's determinant formula.} took a position of a junior assistant in September 1915 at the Astronomical Observatory there, moving from the University of Kazan.  He also submitted his thesis ``Three essays on refraction theory", for which he got habilitation and became a privatdozent at the Yuryev University. It was difficult for him to carry out his observations as planned, because some instruments were evacuated to places further inside the Russian Empire. However, his theoretical work (on orbit determination, involving high-precision solutions to Gauss' equation) was going well and brought him the master's degree in 1917, which in turn led to nomination for a docent, winning the competition for a professor position in 1918 (a vacancy was created by the transfer of Konstantin Dorimedontovich Pokrovskii, 1868-1944, to Perm) and the appointment as the director of the observatory. The university was being transferred to Voronezh  and Banachiewicz got a nomination for a professor's position there, which he did not accept. He was allowed by the German occupying authorities \footnote{The independent Republic of Estonia was declared on February 24, 1918. The Germans withdrew from the territory and handed over control to the Estonian Provisional Government in November 1918.}  to go to Warsaw in September 1918. The end of war and proclamation of independent Poland found him there. From October 1918 to March 1919 he was a deputy professor of geodesy at the Warsaw Polytechnic. In 1919 he took a chair of astronomy at the Jagiellonian University,  which he was offered in May 1918 (\cite{FP},  \cite{koronscy16}).



\section{Poles in other countries}

\textbf{Leon Lichtenstein (1878-1933)} held PhD degrees in engineering and mathematics (from Technische Hochschule Berlin-Char\-lot\-ten\-burg and Fried\-rich-Wilhelm University, respectively) and was active both as a mathematician and an engineer in Berlin. From  1910 (veniam legendi) to 1919 he taught at the Technische Hochschule Berlin-Charlottenburg, lecturing on synthetic and descriptive geometry, graphic static, vector calculus, trigonometric series, integral equations, potential theory and other subjects. At the same time (from  1902 to 1920) he worked for Siemens \& Halske  (later renamed Siemens-Schuckert Werke), becoming a head of the electric laboratory in the factory of electric cables in 1906 and a mathematical expert in 1918.\footnote{Some of Lichtenstein's engineering papers are mentioned here: \textit{High tension cable manufacture, present state and future}, London Electrician, June 2, 1911;
\textit{Testing high-tension cables}, Elek. Zeit., October 8, 1914 .} The electrical industry was important for German economic growth.\footnote{According to \cite{eksteins}, by 1913 the value of German electrical production was twice that of Britain and almost ten times that of France, while Germany's exports in this area were the largest in the world, almost three times those of the United States.} The Siemens company was also active in the arms industry \footnote{This activity resulted in January 1914 in the so-called Siemens Scandal involving bribes for supplying  the Japanese navy.}  and   contributed to the war effort of the German Empire (it developed, among other things, a type of a rotary aircraft engine). The usefulness of his work was probably the reason why Lichtenstein, who  was born in Warsaw and completed one-year ``voluntary" service in the Russian army (in 1897),  was able to obtain German citizenship in the first days of the war. (\cite{przeworska})\\
  
\textbf{Chaim (Herman) M\"untz (1884-1956)}, born in \L \'od\'z, obtained his doctorate in mathematics at the Friedrich-Wilhelm University in Berlin in 1910. He was unsuccessful in getting habilitation and academic position in Germany, so after a period of supporting himself with private lessons he became a teacher of mathematics in a boarding school called the Odenwaldschule near Heppenheim in southern Hessen.\footnote{It was a renowned modern co-educational school, founded and run by the  innovative educator Paul Geheeb (1870-1961).} He was given ample time to work on his mathematical research. In 1915 he left and took a position at another boarding school, also in Hessen (having only Hessian residency but no German citizenship he could not move freely), from which he was dismissed in 1917 as a ``little Polish Jew".   M\"untz was able to write and publish 5 mathematical research papers while teaching.  Also in 1915, he met and befriended the philosopher Martin Buber (1878-1965). He contributed to a journal Der Jude founded and co-edited by Buber (under the pseudonym Herman Glenn).  (\cite{ortizpinkus})\\

\textbf{Mieczys\l aw Biernacki (1891-1959)}: In the years 1909-11 he studied  chemistry at the Jagiellonian University. He was expelled for taking part in students' protests. Then he continued his studies at Sorbonne, switching to mathematics . When the war broke out, he voluntarily enlisted in the French army. He fought at the Western front, suffering gas poisoning and a severe wound. For his service he received the distinction of the Officer's Cross of the Legion of Honour. On June 4, 1917, the president of France issued a decree about forming an independent Polish army in France. Biernacki transferred to the Polish units. He returned to Poland with  the Polish army under the command of General J\'ozef Haller (1873-1960), also known as the ``Blue Army". In 1928 Biernacki obtained a doctorate in Paris under the direction of Paul Montel. (\cite{montel}, \cite{blekitna}, \cite{paryz}). \\

\textbf{Juliusz Pawe\l \ Schauder (1899-1943)}: He graduated from Gymnasium VIII in Lw\'ow in 1917, was drafted into the army  and sent to the Italian front.\footnote{The most important  operations on the Italian front took place in the Isonzo valley. Ultimately the Italians prevailed and the armistice with Austro-Hungary was signed on November 3, 1918.} He was taken a prisoner. While in the camp, he learned about a Polish army being formed in France under the command of Gen. Haller.  On January 24, 1919, he reported to the local recruitment office and was enlisted into a company of ensigns in the rank of corporal. He returned to Poland with Haller's army  and wore his blue uniform long after being discharged, because of material hardship (\cite{derk90}).\\

\textbf{Stefan Straszewicz (1889-1983)} - He taught at the Society for Scientific Courses (fundamental notions of set theory, among other things). In summer 1913 he went to Z\"urich, thanks to the scholarship from the Mia\-now\-ski Fund.\footnote{According to \cite{dziadek}, Stefan Straszewicz's  studies were financed by Zygmunt Straszewicz, his uncle.} He got his PhD at the University of Z\"urich in 1914 under the direction of Ernst Zermelo (1871-1953) on the basis of the thesis  ``Beitr\"age zur Theorie der konvexen Punktmengen" (Research on the theory of convex sets). He continued his research in geometry and topology and translated into Polish the book ``Stetigkeit und irrationale Zahlen" (Continuity and the irrational numbers) by Richard Dedekind (1831-1916). He belonged to the Union of Societies of Polish Youth for Independence (Unia Stowarzysze\'n Polskiej M\l odzie\.zy Niepodleg\l o\'sciowej), commonly called Filarecja. He returned to Poland in 1919. He fought in the Polish-Soviet war, then taught in Warsaw, first at the University and then at the Polytechnic. (\cite{dziadek})\\

\textbf{Acknowledgments:} This article originated from the lecture ``Mathematicians from Polish
territories in WWI"  given by the first author at the conference   ``Mathematics, Mathematicians and World War I",   20
May 2015 - 25 May 2015, Scuola Normale Superiore  Pisa, Italy. It started taking its present shape during the second author's visit to Jagiellonian University in Krak\'ow on her study leave from American Mathematical Society in February-June 2017. The work of both authors is partially supported by the project 18-00449S of Czech Science Foundation. The work of the first author was also  partially supported by the Centre for Innovation and Transfer of Natural
Sciences and Engineering Knowledge, University of Rzesz\'ow.}


\end{document}